# DISPERSION AND SPURIOUS REFLECTION OF VISCOELASTIC WAVE

**José Elias Laier**[*]
[1]Engineering School of São Carlos (EESC), Department of Structural Engineering, University of São Paulo, São Carlos, SP, Brazil

[*]*email:* jelaier@sc.usp.br

**ABSTRACT**

This article investigates the velocity dispersion and the spurious reflection of the viscoelastic wave that occur in the numerical integration of the viscoelastic wave equation. For this purpose, the classic finite element of two nodes, with a consistent and lumped mass model for spatial integration is considered, and the Newmark average acceleration method of the two-step version for integration over time is adopted. The resulting system of the difference equation is then analytically integrated in non-finite terms (numerical solution of waves) using complex notation. The numerical results reveal that, even for a refined mesh, the dispersion and spurious reflections are significant.

## INTRODUCTION

Viscoelastic wave propagation has been a very frequent theme in the technical literature. For example, recently Guo, P.; McMechan, G. A. and Ren, L developed a modified viscoacoustic wave simulation algorithm for modeling the viscoelastic effects of P-waves [1]. Zhong, W and Liu, T suggest a simple triangular mesh grading technique to conduct efficient simulation of viscoelastic wave bidirectional propagation at the irregular interface of large velocity-contrast earth media with rugged topography [2]. Milazzo, M.; Jung, G, S.; Danti, S; et al studied the mechanical resilience and viscoelasticity, wave propagation properties and energy dissipation of collagenous structures [3]. Ferreira, Oliveira and Pena consider general linear damped wave equations with memory [4].

The finite element method solves difference equations (a discrete system) and, just at the limit, the corresponding differential equation. In the case of the wave equation, the solution of the corresponding difference equation presents propagation velocity dispersion [5-7] and spurious reflections for the irregular mesh [8-10]. The numerical integration of the wave equation makes the wave velocity frequency-dependent (dispersive propagation). After a certain time, it is impossible to recognize the shape of a traveling pulse. A spurious wave reflection occurs when a wave passes between two finite elements of different sizes.

Wave propagation in discrete systems has been studied since the 17th century and has been described by Brillouin [11]. A one-dimensional lattice of a point mass interconnected by springs has served as a model for sound wave propagation and the propagation of waves in crystals. Substantial research has been devoted to dispersion and spurious reflections in the numerical integration of non-dispersive wave equations (scalar waves), i.e., waves whose propagation velocity is not frequency-dependent. Belytschko and Mullen [5] present one of the first investigations on the dispersion property of the finite element method solution for one-dimensional waves using linear and quadratic functions for displacement and on the effect of the numerical integration over time on the dispersion. Mullen and Belytschko [6] investigate the dispersion analysis of finite element semidiscretizations of the two-dimensional wave equation considering bilinear quadrilateral elements and triangular elements for various propagation directions. Liu, Sharan, and Yau [7] study the wave motion and its dispersive properties in a finite element model with two-dimensional distortional elements. Bazant [8] discuss the spurious reflection of elastic waves in a non-uniform grid and the effect of consistent and lumped mass models. Bazant and Zekai [9] investigate the spurious reflection of an elastic wave in non-uniform meshes of constant and linear strain finite elements. Jiang and Rogers [9] study the spurious wave reflections at an interface of different physical properties in finite-element wave solutions. Recently, Idesman [10] discussed the optimal reduction in the numerical dispersion for wave propagation problems using two-dimensional isogeometric elements.

Although the study of spurious reflection and dispersion in the numerical integration of scalar waves has been the subject of intensive study, little attention has been paid to the numerical dispersion and spurious reflection of viscoelastic wave that present evanescent propagation. This paper examines the properties of dispersion and spurious reflection in the numerical integration of viscoelastic wave. The spatial integration of the viscoelastic wave equation is performed using the two-node finite element with consistent and lumped mass models, and that over time is performed with the two-step version of the Newmark average acceleration method. The resultant system of the difference equation is then analytically integrated in non-finite terms (numerical wave solution) using complex notation. The numerical results reveal that even for a refined mesh, the dispersion and spurious reflections are significant.

**ONE-DIMENSIONAL VISCOELASTIC WAVE EQUATION**

The classical one dimensional viscoelastic wave equation is expressed by

$$u^{II} + c\dot{u}^{II} - \frac{1}{v_{re}^2}\ddot{u} = 0 \tag{1}$$

where

$$v_r = \sqrt{\frac{E}{\rho}} \tag{2}$$

x and t are the space and time variables, respectively, c is the damping constant, $v_r$ is the reference wave velocity of propagation, E is the modulus of elasticity and $\rho$ is the material density. A Roman numeral as an exponent indicates the degree of the spatial derivative and upper dots represent the time derivative.

Two main techniques are generally used to solve wave equation (1) according to Clough and Penzien [13]. The first technique, which is referred to as integration in finite terms, is applied using normal mode superposition (vibration solution), and the second technique, referred to as integration in non-finite terms, is known as the D´Alembert wave solution [14].

The integration of equation (1) in non-finite terms (wave solution) using complex notation is given by

$$u = Ae^{i(\beta x - \omega t)} \tag{3}$$

where $\beta$ is the wave number, A is the wave amplitude, x is the space variable, t is the time variable, $\omega$ is the angular wave frequency and i is the complex a unit. Substituting (3) into (1) results in the following:

$$\beta^2 - \omega c \beta^2 i - \frac{\omega^2}{v_r^2} = 0 \tag{4}$$

The roots of (4) are given by the following:

$$\beta = \pm \frac{\omega}{v_r}\sqrt{\frac{1}{1-\omega ci}} = \pm \frac{\omega}{v_r}\sqrt{\frac{1+\omega ci}{1+\omega^2 c^2}} = \pm \frac{\omega}{v_r}\sqrt{\overline{A} + \overline{B}i} \tag{5}$$

where

$$\begin{aligned}\overline{A} &= \frac{1}{1+\omega^2 c^2} \\ \overline{B} &= \frac{\omega c}{1+\omega^2 c^2}\end{aligned} \tag{6}$$

Equation (5) in the polar form can be expressed as follows:

$$\beta = \frac{\omega}{v_r}\sqrt{\rho_c e^{\theta i}} \tag{7}$$

where

$$\rho_c = \frac{1}{\sqrt{1+\omega^2 c^2}}$$
$$\theta = \arc\left(\cos = \frac{\sqrt{1+\omega^2 c^2}}{1+\omega^2 c^2}\right) \tag{8}$$

$\rho_c$ is the modulus of the complex wave number and $\theta$ is the phase. Note that in equation (3), without any loss of generality, just the wave propagation in the positive direction is considered. Equation (7) can be written as

$$\beta = \pm\frac{\omega}{v_r}\sqrt[4]{\frac{1}{1+\omega^2 c^2}}\left(\cos\left(\frac{\theta}{2}\right)+i\operatorname{sen}\left(\frac{\theta}{2}\right)\right) = A^* + B^* i \tag{9}$$

where

$$A^* = \frac{\omega}{v_r}\sqrt[4]{\frac{1}{1+\omega^2 c^2}}\cos\left(\frac{\theta}{2}\right)$$
$$B^* = \frac{\omega}{v_r}\sqrt[4]{\frac{1}{1+\omega^2 c^2}}\operatorname{sen}\left(\frac{\theta}{2}\right) \tag{10}$$

Finally, taking in to account equation (9) equation (3) can be expressed by

$$u = Ae^{i\left[\left(A^*+B^*i\right)x-\omega t\right]} = Ae^{-B^* x}e^{i\left(A^* x-\omega t\right)} \tag{11}$$

corresponding to a wave with propagation velocity of

$$v = \frac{\omega}{A^*} \tag{12}$$

which means that the velocity of propagation is frequency-dependent (dispersive wave) and decays exponentially over of the space variable (evanescent propagation). $A^*$ and $B^*$ are the wave number and wave damping, respectively.

**FINITE ELEMENT FORMULATION**

The two-node element formulation, as illustrated in Figure 1, can be summarized by the following matrix equation [15]:

$$\frac{E}{\ell}\begin{bmatrix}1 & -1\\-1 & 1\end{bmatrix}\begin{Bmatrix}u_j\\u_{j+1}\end{Bmatrix} + c\frac{E}{\ell}\begin{bmatrix}1 & -1\\-1 & 1\end{bmatrix}\begin{Bmatrix}\dot{u}_j\\\dot{u}_{j+1}\end{Bmatrix} + \rho\ell\begin{bmatrix}M_1 & M_2\\M_2 & M_1\end{bmatrix}\begin{Bmatrix}\ddot{u}_j\\\ddot{u}_{j+1}\end{Bmatrix}$$
$$-\begin{Bmatrix}-N_j\\N_{j+1}\end{Bmatrix} = \begin{Bmatrix}0\\0\end{Bmatrix} \qquad (13)$$

where $\ell$ is the element length, $M_1 = 1/3$ and $M_2 = 1/6$ for a consistent mass model, $M_1 = 1/2$ and $M_2 = 0$ for a classical lumped mass matrix, and $N_j$ and $N_{j+1}$ are normal forces at nodes j and j+1, respectively, as in Figure 1.

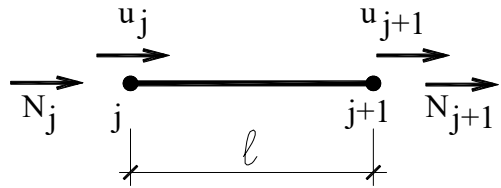

Figure 1 Two node Finite Element

By considering equation (4), the equilibrium of a generic node j, which consists of the finite element version of wave equation (1), can be expressed by

$$\{-u_{j-1} + 2u_j - u_{j+1}\} + \frac{c}{\Delta t}\{-\dot{u}_{j-1} + 2\dot{u}_j - \dot{u}_{j+1}\} + \frac{\rho\ell^2}{E}\{M_2\ddot{u}_{j-1} + 2M_1\ddot{u}_j + M_2\ddot{u}_{j+1}\} = 0 \qquad (14)$$

as shown in Figure 2.

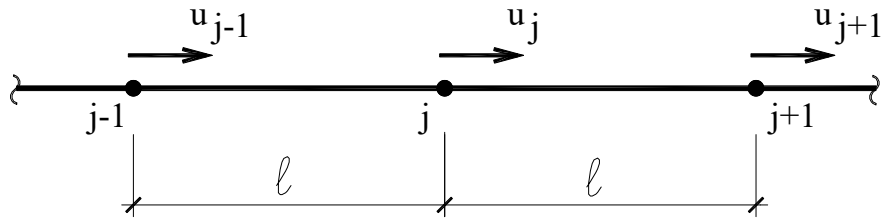

Figure 2 Finite Element Equilibrium

Note that equation (14) is a differential difference equation (differential in time and difference in space) in which the space variable x is replaced by the discrete variable $j\ell$ where j=0,1… and $\ell$ is now considered the space increment.

**TIME INTEGRATION ALGORITHM**

For time integration, the two-step version of Newmark's average acceleration method is considered, and thus, it is expressed as

$$\dot{u}_{k+1} - 2\dot{u}_k + \dot{u}_{k-1} + \frac{\Delta t}{2}\left(\ddot{u}_{k-1} - \ddot{u}_{k+1}\right) = 0$$

$$u_{k+1} - 2u_k + u_{k-1} - \frac{\Delta t^2}{4}\left(\ddot{u}_{k-1} + 2\ddot{u}_k + \ddot{u}_{k+1}\right) = 0 \tag{15}$$

Note that equation (15) is a difference equation in time involving three unknown functions ($u$, $\dot{u}$ and $\ddot{u}$) in which the time variable t is replaced by the discrete variable k$\Delta$t where k=0,1… and $\Delta$t is the time increment. It is necessary to consider two independent Hermitian operators because the structural dynamic equations are second-order in time.

Equation (14) and equation (15) form a system of three difference equations, now in discrete space and time variables, involving those three discrete unknowns functions, and whose solutions are given by

$$u_{j,k} = A e^{i(\beta_n j \ell - \omega k \Delta t)}$$
$$\dot{u}_{j,k} = B e^{i(\beta_n j \ell - \omega k \Delta t)} \tag{16}$$
$$\ddot{u}_{j,k} = C e^{i(\beta_n j \ell - \omega k \Delta t)}$$

where $\beta_n$ is the numerical wave number and $A$, $B$ and $C$ are the amplitude of the displacement, velocity and acceleration, respectively. Note that equation (16) permits us to form, for example, the following relationships that stem from Euler's formulas

$$u_{j-1,k} + u_{j+1,k} = A e^{i(\beta_n j \ell - \omega k \Delta t)}\left[2\cos(\beta_n \ell)\right]$$
$$u_{j,k-1} - u_{j,k+1} = A e^{i(\beta_n j \ell - \omega k \Delta t)}\left[2i\,\mathrm{sen}(\omega \Delta t)\right] \tag{17}$$
$$u_{j,k-1} + u_{j,k+1} = A e^{i(\beta_n j \ell - \omega k \Delta t)}\left[2\cos(\omega \Delta t)\right]$$

**VELOCITY DISPERSION ANALYSIS**

Taking into account the expression in (16) and (17), equations (14) and (15) are rewritten in matrix notation as follows:

$$\begin{bmatrix} 1-\cos(\beta_n \ell) & \psi_2(1-\cos(\beta_n \ell)) & \psi_1(M_2 \cos\beta_n \ell + M_1) \\ 0 & \cos(\omega\Delta t)-1 & \dfrac{i}{2}\mathrm{sen}(\omega\Delta t) \\ \cos(\omega\Delta t)-1 & 0 & -\dfrac{\cos(\omega\Delta t)+1}{4} \end{bmatrix} \begin{Bmatrix} A \\ B\Delta t \\ C\Delta^2 t \end{Bmatrix} = \begin{Bmatrix} 0 \\ 0 \\ 0 \end{Bmatrix} \tag{18}$$

where

$$\Psi_1 = \frac{\rho \ell^2}{E \Delta t^2} = \left(\frac{a}{b}\right)^2$$

$$\Psi_2 = \frac{c}{\Delta t} = \frac{ac}{T} \tag{19}$$

and

$$a = \frac{T}{\Delta t}$$

$$b = \frac{\lambda}{\ell} \tag{20}$$

In addition a is the mesh parameter in time, b is the mesh parameter in space, T is the wave period and $\lambda = v_r T$ is the reference wave length ($v_r$ is the reference propagation velocity).

The eigen-value of the eigen-value problem expressed by equation (18) is given by

$$\cos(\beta_n \ell) = \frac{D_1 + F_1 i}{D_2 + F_1 i} \tag{21}$$

where

$$D_1 = -\left[\frac{\cos(\omega \Delta t) + 1}{4} + \psi_1 M_1 (\cos(\omega \Delta t) - 1)\right]$$

$$D_2 = \left[-\frac{\cos(\omega \Delta t) + 1}{4} + \psi_1 M_2 (\cos(\omega \Delta t) - 1)\right] \tag{22}$$

$$F_1 = \frac{\psi_2 \mathrm{sen}(\omega \Delta t)}{2}$$

and

$$\cos(\beta_n \ell) = D + F i \tag{23}$$

where

$$D = \left[D_1 D_2 + \frac{\psi_2^2 \mathrm{sen}^2(\omega \Delta t)}{4}\right] / \left[D_2^2 + \frac{\psi_2^2 \mathrm{sen}^2(\omega \Delta t)}{4}\right]$$

$$F = \left[\frac{\psi_1 \psi_2 (\cos(\omega \Delta t) - 1)(M_1 + M_2) \mathrm{sen}(\omega \Delta t)}{2}\right] / \left[D_2^2 + \frac{\psi_2^2 \mathrm{sen}^2(\omega \Delta t)}{4}\right] \tag{24}$$

The solution of the transcendental equation (23) is expressed by the following:

$$\beta_n \ell = d + hi \qquad (25)$$

where

$$\cos(d) = \sqrt{\frac{1+D^2+F^2 - \sqrt{(1+D^2+F^2)-4D^2}}{2}}$$
$$\cosh(h) = \sqrt{\frac{1+D^2+F^2 + \sqrt{(1+D^2+F^2)-4D^2}}{2}} \qquad (26)$$

and

$$d = \text{arc}(\cos) = \sqrt{\frac{1+D^2+F^2 - \sqrt{(1+D^2+F^2)-4D^2}}{2}}$$
$$h = \ell n\left(\cosh(h) + \sqrt{\cosh^2(h)-1}\right) \qquad (27)$$

Note that the numerical wave number $\beta_n = d/\ell$ and numerical wave damping $h/\ell$ are frequency and damping dependent (dispersive wave).

Finally, the eigen-vector of equation (18) assuming a unitary amplitude for the wave displacement can be written as follows:

$$\begin{Bmatrix} A \\ B\Delta t \\ C\Delta^2 t \end{Bmatrix} = \begin{Bmatrix} u_{j,k} \\ \Delta t \dot{u}_{j,k} \\ \Delta t^2 \ddot{u}_{j,k} \end{Bmatrix} = \begin{Bmatrix} 1 \\ -2i\text{sen}(\omega\Delta t)/(\cos(\omega\Delta t)+1) \\ 4(\cos(\omega\Delta t)-1)/(\cos(\omega\Delta t)+1) \end{Bmatrix} \qquad (28)$$

**VISCOELASTIC WAVE DISPERSION AND NUMERICAL DAMPING RESULTS**

Tables 1 and 2 present the results for the percentage of the relative numerical velocity error and the relative numerical damping error. A one-dimensional viscoelastic homogeneous medium subdivided by two-node finite element is adopted. Consistent and lumped mass models are considered. The lumped mass model results are given in brackets. The percentage of the relative numerical wave velocity error, relative numerical damping error and physical damping are expressed, respectively, as

$$100\frac{v - v_n}{v} = 100\frac{d/\ell - A^*}{d/\ell}$$

$$100\frac{B^* - h/\ell}{B^*} \qquad (29)$$

$$\gamma = \frac{\pi \Psi_2}{a}$$

where $v_n$ is the numerical wave velocity. Table 1 includes equal mesh parameters in space and time (unity Courant number $b/a = 1$) and three different physical damping values. Table 2 present the results for different values of the mesh parameters in space and time (different Courant numbers).

Table 1 Percentage of the relative numerical wave velocity and damping errors

| MESH PARAMETERS | PHYSICAL DAMPING | RELATIVE NUMERICAL DAMPING ERROR % | RELATIVE NUMERICAL WAVE VELOCITY ERROR% |
|---|---|---|---|
| a=b=100 | $\gamma$=0.1 | 0.01645 (0.04747) | 0.01646 (0.1121) |
| | $\gamma$=0.01 | 0.01645 (0.04934) | 0.01646 (0.1153) |
| | $\gamma$=0.001 | 0.01645 (0.04936) | 0.01669 (0.1158) |
| a=b=50 | $\gamma$=0.1 | 0.06582 (01900) | 0.06592 (0.4500) |
| | $\gamma$=0.01 | 0.06582 (0.1975) | 0.06593 (0.4627) |
| | $\gamma$=0.001 | 0.06582 (0.1976) | 0.06592 (0.4628) |
| a=b=25 | $\gamma$=0.1 | 0.2635 (0.7620) | 0.2652 (1.825) |
| | $\gamma$=0.01 | 0.2635 (0.7924) | 0.2654 (1.878) |
| | $\gamma$=0.001 | 0.2635 (0.7927) | 0.2654 (1.878) |
| a=b=20 | $\gamma$=0.1 | 0.4119 (1.193) | 0.4161 (2.883) |
| | $\gamma$=0.01 | 0.4121 (1.241) | 0.4166 (2.967) |
| | $\gamma$=0.001 | 0.4121 (1.241) | 0.4166 (2.968) |
| a=b=10 | $\gamma$=0.1 | 1.657 (4.857) | 1.731 (12.68) |
| | $\gamma$=0.01 | 1.660 (5.065) | 1.740 (13.10) |
| | $\gamma$=0.001 | 1.660 (5.967) | 1.740 (13.10) |

Although the viscoelastic wave is dispersive (velocity of propagation is frequency-dependent), the relative numerical wave velocity error and the relative numerical damping error for unity Courant number ($a = b$) practically do not vary with the wave frequency. Table 1 shows that the effect of physical damping on the results can be negligible. Note that the results of the lumped mass model for the relative numerical damping error are something like three times greater than those of the consistent mass model, but the relative numerical wave velocity errors of the lumped mass model are seven times greater than those of the consistent mass model. Furthermore, Table 1 shows that, in the case of a coarse mesh ($a = b = 10$), the relative numerical wave velocity error and the relative numerical damping error for the lumped mass model are quite high ($\approx 10\%$).

Table 2 Relative numerical wave velocity and damping errors in percentages

| MASH PARAMETERS | RELATIVE NUMERICAL DAMPING ERROR % | RELATIVE NUMERICAL WAVE VELOCITY ERROR % |
|---|---|---|
| a=100   b=50 | -0.030002 (0.09405) | -0.1265 (0.2561) |
| a=50   b=100 | 0.1123 (0.1434) | 0.2094 (0.3054) |
| a=50   b=25 | -0.1195 (0.3772) | -0.5021 (1.034) |
| a=25   b=50 | 0.4495 (0.574) | 0.8419 (1.232) |
| a=20   b=10 | -0.7238 (2.403) | -2.971 (6.944) |
| a=10   b=20 | 2.818 (3.617) | 5.453 (8.176) |
| a=10   b=5 | -2.566 (10.42) | -9.724 (38.43) |
| a=5   b=10 | 11.39 (14.93) | 24.98 (41.88) |

The results listed in Table 2 show for a fine mesh ($a, b > 50$) that the relative numerical damping error and relative numerical wave velocity error are less than 1% for the consistent mass model and for the lumped mass model. However, for a coarse mesh ($a, b < 20$) the results for the relative numerical damping error are high ($\approx 10\%$) and the results for the relative numerical wave velocity error are quite high ($\approx 40\%$).

**SPURIOUS WAVE REFLECTIONS**

Consider a non-uniform finite element mesh configuration shown in Figure 3. By considering an incident wave with unitary amplitude traveling from left to right and arriving at the interface, a spurious reflected wave and a transmitted wave are generated.

The incident wave in complex notation is thus expressed as

$$u_{in} = \exp[i(\beta_{in} j\ell - \omega k \Delta t)] \tag{30}$$

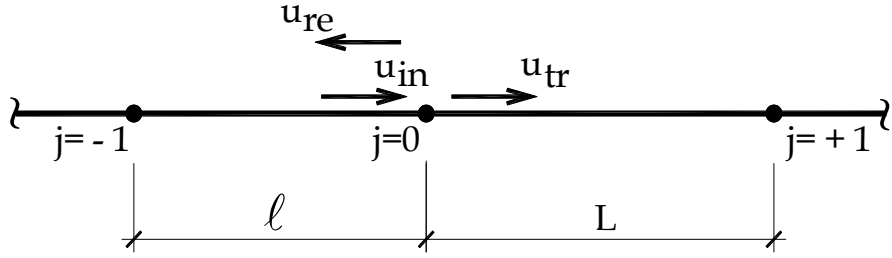

Figure 3 Non-uniform Finite Element Mesh

where $\beta_{in}$ is the incident numerical wave number. The spurious reflected wave is given by

$$u_{re} = A_{re} \exp[i(\beta_{in} j\ell - \omega k \Delta t)] \quad (31)$$

where $A_{re}$ is the amplitude of the spurious reflected wave. The transmitted wave is expressed as

$$u_{tr} = A_{tr} \exp[i(\beta_{in} j\ell - \omega k \Delta t)] \quad (32)$$

where $A_{tr}$ is the amplitude of the transmitted wave and $\beta_{tr}$ is the numerical wave number of the right side mesh (Figure 3).

To calculate the spurious reflection amplitude $A_{re}$ and transmitted wave amplitude $A_{tr}$, one should consider the displacement compatibility at the interface (node j), which is expressed by

$$1 + A_{re} = A_{tr} \quad (33)$$

The equilibrium at node j is written as

$$\{-1 \quad 1\}\begin{Bmatrix} u_{-1} \\ u_0 \end{Bmatrix} + \left\{\frac{1}{\alpha} \quad -\frac{1}{\alpha}\right\}\begin{Bmatrix} u_0 \\ u_{+1} \end{Bmatrix} + \psi_2\left[\{-1 \quad 1\}\begin{Bmatrix} \Delta t \dot{u}_{-1} \\ \Delta t \dot{u}_0 \end{Bmatrix} + \left\{\frac{1}{\alpha} \quad -\frac{1}{\alpha}\right\}\begin{Bmatrix} \Delta t \dot{u}_0 \\ \Delta t \dot{u}_{+1} \end{Bmatrix}\right] +$$
$$\psi_1\left[\{M_1 \quad M_2\}\begin{Bmatrix} \Delta t^2 \ddot{u}_{-1} \\ \Delta t^2 \ddot{u}_0 \end{Bmatrix} + \alpha\{M_2 \quad M_1\}\begin{Bmatrix} \Delta t^2 \ddot{u}_0 \\ \Delta t^2 \ddot{u}_1 \end{Bmatrix}\right] = 0 \quad (34)$$

where

$$\alpha = \frac{L}{\ell}$$
$$u_{-1} = [\exp(-i\beta_{in}\ell) + A_{re}\exp(i\beta_{in}\ell)]\exp(-i\omega k \Delta t)$$
$$u_0 = (1 + A_{re})\exp(-i\omega k \Delta t) = A_{tr}\exp(-i\omega k \Delta t) \quad (35)$$
$$u_{+1} = A_{tr}\exp(i\beta_{tr}L)\exp(-i\omega k \Delta t) = (1 + A_{re})\exp(i\beta_{tr}L)\exp(-i\omega k \Delta t)$$

Substituting (34) in (33), and taking into account (28), after algebraic manipulation, we get the following:

$$A_{re} = -\frac{(\bar{A}+A^*i)e^{-i\beta\ell}+(\bar{B}+B^*i)+(\bar{C}+C^*i)e^{i\beta_{tr}L}}{(\bar{A}+A^*i)e^{i\beta\ell}+(\bar{B}+B^*i)+(\bar{C}+C^*i)e^{i\beta_{tr}L}} \qquad (36)$$

where

$$\bar{A} = -1+4\psi_1 \frac{\cos(\omega\Delta t)-1}{1+\cos(\omega\Delta t)}M_2$$

$$A^* = 2\psi_2 \frac{\operatorname{sen}(\omega\Delta t)}{1+\cos(\omega\Delta t)}$$

$$\bar{B} = 1+\frac{1}{\alpha}+4\psi_1 M_1 \frac{\cos(\omega\Delta t)-1}{1+\cos(\omega\Delta t)}(1+1/\alpha)$$

$$B^* = -2\psi_2 \frac{\operatorname{sen}(\omega\Delta t)}{1+\cos(\omega\Delta t)}(1+1/\alpha) \qquad (37)$$

$$\bar{C} = -\frac{1}{\alpha}+4\psi_1\alpha \frac{\cos(\omega\Delta t)-1}{1+\cos(\omega\Delta t)}M_2$$

$$C^* = \frac{2}{\alpha}\psi_2 \frac{\operatorname{sen}(\omega\Delta t)}{1+\cos(\omega\Delta t)}$$

On the other hand, expression (23) allows one to write

$$\cos(\beta_{in}\ell) = D - Fi$$
$$\cos(\beta_{tr}L) = \bar{D} - \bar{F}i \qquad (38)$$

and

$$\operatorname{sen}(\beta_{in}\ell) = G + Hi$$
$$\operatorname{sem}(\beta_{tr}\ell) = \bar{G} + \bar{H}i \qquad (39)$$

where

$$H = \frac{-(1-D^2+F^2)+\sqrt{(1-D^2+F^2)^2+4D^2F^2}}{2}$$

$$G = \frac{DF}{H}$$

$$\bar{H} = \frac{-(1-\bar{D}^2+\bar{F}^2)+\sqrt{(1-\bar{D}^2+\bar{F}^2)^2+4\bar{D}^2\bar{F}^2}}{2} \qquad (40)$$

$$\bar{G} = \frac{\bar{D}\bar{F}}{\bar{H}}$$

After algebraic manipulations, expression (35) becomes the following:

$$A_{re} = -\frac{d_1d_3 + d_2d_4 + i(d_2d_3 - d_1d_4)}{d_3^2 + d_4^2} \quad (41)$$

where

$$\begin{aligned}
d_1 &= \overline{A}D + A^*G + \overline{B} + \overline{C}\overline{D} - C^*\overline{G} + A^*F + \overline{A}H + C^*\overline{F} - \overline{C}\overline{H} \\
d_2 &= -\overline{A}F + A^*H - \overline{C}\overline{F} - C^*\overline{H} + A^*D - \overline{A}G + B^* + C^*\overline{D} + \overline{C}\overline{G} \\
d_3 &= \overline{A}D - A^*G + \overline{B} + \overline{C}\overline{D} - C^*\overline{G} + A^*F - \overline{A}H + C^*\overline{F} - \overline{C}\overline{H} \\
d_4 &= -\overline{A}F - A^*H - \overline{C}\overline{F} - C^*\overline{H} + A^*D + \overline{A}G + B^* + C^*\overline{D} + \overline{C}\overline{G}
\end{aligned} \quad (42)$$

The spurious reflection module is then the following:

$$|A_{re}| = \frac{\sqrt{(d_1d_3 + d_2d_4)^2 + (d_2d_3 - d_1d_4)^2}}{d_3^2 + d_4^2} \quad (43)$$

which is what matters in practical terms.

Table 3 shows the results of the reflected spurious wave for various mesh parameters (with different values for the alpha ratio). Note that the results for both the consistent and lumped mass models are quite similar. However, the amplitude of the spurious wave has values above 1%, even for a refined mesh ($a = b = 100$) with an alpha ratio slightly away from the unit. On the other hand, even for moderately refined meshes, the amplitude of the reflected spurious wave is significant. The amplitude of the spurious reflection has practically the same magnitude as the incident wave for a coarse mesh ($a = 10$, $b = 5$).

**CONCLUSION**

The wave velocity dispersion and spurious wave reflections at the interfaces of finite elements of different lengths (non-uniform mesh) are examined here. The numerical results reveal that the consistent mass model is superior to the lumped mass model as expected. The consistent mass model allows better resolution for short wavelengths and more frequency representation, as is well known.

The relative error of the numerical damping and the relative and numerical wave velocity are less than 1% for mesh parameters greater than 20 in the case of the consistent mass model and greater than 50 in the case of the lumped mass model (meshes with a unit Courant number). For meshes with a non-unity Courant number, the mesh parameters $a$ and $b$ must be greater than 50.

The amplitude of the reflected spurious wave is less than 1% of the incident wave amplitude for the mesh parameter a = 100 and an alpha ratio in the range of 0.9 to 1.1. For other values of the mesh parameters and ranges of alpha ratio, the reflected spurious wave can reach as much as 70% of the incident wave.

Table 3 Relative spurious wave amplitude

| MESH PARAMETER | α | MESH PARAMETER | SPURIOUS REFLECTIONS % |
|---|---|---|---|
| a=100 | 0.5 | 200 | 1.845 (1.845) |
| | 0.7 | 142.9 | 1.397 (1.397) |
| | 0.9 | 111.1 | 0.5740 (0.5753) |
| | 1.1 | 90.91 | 0.6966 (0.6937) |
| | 1.5 | 66.67 | 4.868 (4.860) |
| | 2.0 | 50 | 13.67 (13.65) |
| a=50 | 0.5 | 100 | 3.704 (3.705) |
| | 0.7 | 71.35 | 2.798 (2.803) |
| | 0.9 | 55,56 | 1.151 (1.163) |
| | 1.1 | 45.45 | 1.392 (1.371) |
| | 1.5 | 33.33 | 9.605 (9.548) |
| | 2.0 | 25 | 25.86 (25.72) |
| a=10 | 0.5 | 20 | 19.26 (20.10) |
| | 0.7 | 14.29 | 14.39 (15.55) |
| | 0.9 | 11.11 | 5.738 (7.542) |
| | 1.1 | 9.091 | 6.678 (4.080) |
| | 1.5 | 6.667 | 37.51 (33.08) |
| | 2.0 | 5.000 | 66.34 (61.66) |

**ACKNOWLEDGEMENTS**

The authors gratefully acknowledge the funding contributions from the Brazilian National Council for Technological and Scientific Development (CNPq).